\documentclass{article}
\usepackage{amsmath}
\usepackage{amscd}
\usepackage{latexsym}
\usepackage{amssymb}

\title{Self-dual and quasi self-dual algebras}
\author{M. GERSTENHABER}

\begin{document}
\maketitle
\newtheorem{theorem}{Theorem}
\newtheorem{corollary}{Corollary}
\newtheorem{lemma}{Lemma}
%{\theorembodyfont{\rmfamily} \newtheorem{Rem}{Remark}}
\renewcommand{\abstractname}{}
\newcommand{\C}{\ensuremath{\mathbb{C}}}
\newcommand{\Z}{\ensuremath{\mathbb{Z}}}
\newcommand{\I}{\ensuremath{\mathcal{P}}}
\newcommand{\pr}{\ensuremath{\preceq}}
\newcommand{\op}{\ensuremath{\mathrm{op}}}
\newcommand{\g}{\frak{g}}
\newcommand{\G}{\frak{G}}
\vspace{-7mm}
\date{}
{\noindent ${}^1$\textit{Department of Mathematics, University of
Pennsylvania, Philadelphia, PA 19104-6395, \linebreak[0]
U.S.A.\linebreak[0] email:mgersten@math.upenn.edu}

\begin{abstract}\noindent  \emph{Self-dual}
algebras are ones with an  $A$
bimodule isomorphism $A \to A^{\vee\op}$, where $A^{\vee} =
\operatorname{Hom}_k(A,k)$ and $A^{\vee\op}$ is the the same underlying 
$k$-module as $A^{\vee}$ but with left and right operations by $A$ interchanged.
These are in particular \emph{quasi self-dual} algebras,  i.e., ones with an 
isomorphism $H^*(A,A) \cong H^*(A, A^{\vee\,\op})$. For all
such algebras $H^*(A,A)$ is a contravariant functor of $A$.  Finite dimensional
associative self-dual algebras over a field are identical with
symmetric Frobenius algebras. (The monoidal category of commutative 
Frobenius algebras is known to be equivalent to that of  1+1
dimensional topological quantum field theories.) All finite poset algebras
are quasi self-dual.
\end{abstract}

The cohomology of an algebra $A$ with coefficients in itself generally has no functorial properties, but there is an important class of algebras for which $H^*(A,A)$ is  a contravariant functor of $A$. These are the \emph{quasi self-dual algebras}, ones for which there is an isomorphism  $H^*(A,A) \cong H^*(A, A^{\vee\,\op})$,  where $A^{\vee} =\operatorname{Hom}_k(A,k)$ is naturally an $A^{\op}$-bimodule and $A^{\vee\op}$ is the the same underlying $k$-module as $A^{\vee}$ but as  an $A$-bimodule has left and right operations by $A$ interchanged.  \emph{Self-dual algebras}, ones where there is a already an $A$-bimodule isomorphism $A\cong A^{\vee\op}$, are in particular quasi self-dual. (The isomorphisms need not be unique but often there is a natural choice.)  For finite dimensional associative algebras over a field, the category of self-dual algebras is identical with that of symmetric Frobenius algebras, as we show.  Abrams \cite{Abrams:Frobenius}  shows that the  category of commutative Frobenius algebras is equivalent, as a monoidal category, to that of 1+1 dimensional topological quantum field theories (TQFTs), a previous `folk theorem'.  Other papers useful for understanding the relation between Frobenius algebras and TQFTs include those of Sawin, \cite{Sawin:DirectSum}, 
\cite{Sawin:LinksQGsTQFts}; 
the former contains some of the results reproven more directly in \cite{Abrams:Frobenius}. Since $H^*(A,A)$ governs the deformations of $A$, for self-dual and quasi self-dual algebras,  a morphism between algebras may give some relation between their deformation theories.  The connection to TQFTs may be related to the functoriality of $H^*(A,A)$ for such algebras.

For categories of algebras where it is meaningful to consider the
cohomology $H^*(A,k)$ of $A$ with coefficients in $k$ as a trivial
module this cohomology is naturally contravariant and frequently does have a geometric
interpretation. For example, if $G$ is
a finite group operating freely on a contractible space $\mathcal{S}$ 
then the Hochschild cohomology $H^*(kG,k)$
is naturally isomorphic with the cohomology of the quotient space
$\mathcal{S}/G$ with coefficients in $k$. However, $H^*(kG, kG)$ does not have such a simple
geometric interpretation, and while $H^*(kG, kG)$ and $H^*(kG,k)$
are both in a natural way rings, there does not seem to be any
simple way to deduce the structure of the former from the latter.
Nevertheless finite group rings are self-dual as we show, so, in particular, any morphism $G\to G'$
of finite groups induces a morphism $H^*(kG', kG') \to H^*(kG, kG)$.  This 
functoriality raises the question of whether there is  a 
 geometric interpretation of any self-dual or quasi self-dual algebra.
Another basic question not considered here is, When does a morphism
$B\to A$ of self-dual or quasi self-dual algebras induce a morphism $H^*(A,A)\to
H^*(B,B)$ not merely of $k$-modules but one which preserves
additional structures which $H^*(A,A)$ may possess such as, in the
associative case, the Gerstenhaber algebra structure.

We consider here
only algebras $A$ which are finite free modules over their coefficient rings $k$, and
likewise for $A$-bimodules, but some of the most interesting cases
are likely to be infinite dimensional, requiring topological
considerations. On the other hand, the ideas apply both to associative algebras and Lie algebras (and possibly to ones over other operads).

The following problem motivated this paper. To every finite poset
(partially ordered set) $\I$, which we may view as a small category,
one can associate a $k$-algebra $A$ such that $H^*(A,A)$ is
canonically isomorphic with the simplicial cohomology of the nerve of
$\I$ with coefficients in $k$; the proof will be revisited below. If
we start with a simplicial complex $\Sigma$ then its faces, ordered
by setting $\sigma \preceq \tau$ whenever $\sigma$ is a face of
$\tau$, form a poset whose nerve is the barycentric subdivision of
$\Sigma$ and which therefore has the same simplicial cohomology.
(Amongst the faces of a simplex here one includes the simplex
itself.) To every finite simplicial complex $\Sigma$ one can
therefore associate a $k$-algebra $A =A(\Sigma)$ such that $H^*(A,
A) \cong H^*(\Sigma, k)$, where the left side is Hochschild
cohomology and the right simplicial. (Simplicial cohomology is in
fact a special case of Hochschild cohomology independent of any
finiteness assumption, cf. \cite{GS:Simplicial}.) Using the
isomorphism with simplicial cohomology one sees that if $A$ and $B$
are poset algebras then a poset morphism $B\to A$ in fact induces a
morphism $H^*(A,A) \to H^*(B,B)$ just because it induces a
simplicial morphism of the associated simplicial complexes. The
problem was to exhibit a purely algebraic reason for this
functoriality. We show (by a simple extension of the proof of the
preceding result) that poset algebras are in fact quasi self-dual.
As a consequence it is not necessary to restrict to algebra
morphisms induced by morphisms of the underlying posets. There are
other morphisms between poset algebras but they are probably very
restricted in nature.

\section{Opposite algebras and modules.} The \emph{opposite} of an
associative $k$ algebra $A$, denoted $A^{\op}$, is the same
underlying $k$-module but with reversed multiplication: $a\circ b
\in A^{\op}$ is defined to be the element $ba\in A$. (We will
generally use ``$\circ$'' to denote an opposite multiplication.)
Left modules over $A$ are right modules over $A^{\op}$ and an
$A$-bimodule may be viewed as a left module over the universal
enveloping algebra $A^{\mathrm{e}} = A \otimes A^{\op}$. An algebra
morphism of the form $\phi:B\to A^{\op}$ is sometimes called an
antimorphism; it is a $k$-module mapping such that $\phi(ab) = \phi
b\cdot\phi a$ (and is the same thing as a morphism from $B^{\op}$ to
$A$). If $\phi:B\to A$ is an antimorphism and $M$ an $A$-bimodule
then the underlying $k$-module of $M$ becomes a $B^{\op}$-bimodule
by setting $b\circ m \circ b' = \phi b' \cdot m \cdot\phi b$. The
universal enveloping algebra $A^{\mathrm{e}}$ has an
antiautomorphism interchanging the tensor factors. If $G$ is a
group, then its group algebra $kG$ has an antiautomorphism $\sigma$
sending $g\in G$ to $g^{-1}$, extended linearly. The infinitesimal
version of this is that every Lie algebra $\g$ has an
antiautomorphism sending $a\in \g$ to $-a$.

The \emph{opposite of an }$A$-\emph{bimodule} $M$ is the 
$A^{\op}$-bimodule $M^{\op}$ which is the same underlying $k$-module but where
we set $a\circ m\circ b\in M^{\op}$ equal to $b\,m\,a \in M$. Note
that $M^{\op}$ is an $A^{\op}$-bimodule. Although a commutative
algebra $A$ is identical with its opposite, a bimodule $M$ over $A$
may be distinct from its opposite since the left and right
operations of $A$ may be different, as in $A^e$ considered as an
$A$-bimodule. If they are identical then $M$ is called
\emph{symmetric}. Left modules over a commutative algebra may be
viewed as symmetric bimodules.

\begin{theorem} Let $M$ be an $A$-bimodule and $F \in C^n(A,M)$ be a
Hochschild $n$-cochain. Define $F^{\op}\in C^n(A^{\op},M^{\op})$ by
$F^{\op}(a_1,\dots,a_n) = F(a_n,\dots,a_1)$. The map $C^*(A,M) \to
C^*(A^{\op}, M^{\op})$ defined by $F \mapsto
(-1)^{\lfloor\frac{n+1}{2}\rfloor}F^{\op}$ is a cochain isomorphism
inducing an isomorphism of cohomology $H^*(A,M) \cong H^*(A^{\op},
M^{\op})$.
\end{theorem}
\textsc{Proof}. One has
\begin{multline*}
\delta(F^{\op})(a_1,\dots, a_{n+1})=\\
a_1\circ F^{\op}(a_2,\dots, a_{n+1}) -F^{\op}(a_1\circ a_2, \dots,
a_{n+1}) + \dots +(-1)^{n+1}F^{\op}(a_1, \dots, a_n)\circ a_{n+1}\\
=F(a_{n+1},\dots, a_2)a_1 -F(a_{n+1},\dots,
a_3,a_2 a_1) + \dots +(-1)^{n+1}a_{n+1}F(a_n, \dots, a_1) \\
=(-1)^{n+1}(\delta F)^{\op}(a_1,\dots, a_{n+1}).
\end{multline*}
The introduction of the sign $(-1)^{\lfloor\frac{n+1}{2}\rfloor}$
just corrects for the sign $(-1)^{n+1}$ in the foregoing. $\square$
\medskip

Similarly, if $\phi: B\to A$ is an antimorphism and if $M$ is an $A$
bimodule then the $k$-module morphism $\phi^*:C^*(A,M) \to
C^*(B^{\op},M)$ sending $F\in C^n(A,M)$ to $\phi^n F \in
C^n(B^{\op}, M)$  defined by
$$(\phi^n F) (b_1, \dots, b_n) = (-1)^{\lfloor
\frac{n+1}{2}\rfloor}\phi(F(\phi b_n, \dots, \phi b_1)$$ is a
cochain morphism. Likewise, if $M, N$ are $A$ bimodules and $T:M\to
N^{\op}$ is an antimorphism then the $k$-module morphism
$T^*:C^*(A,M) \to C^*(A^{\op},N)$ sending $F\in C^n(A,M)$ to $T^nF
\in C^n(A^{\op},N)$ defined by
$$(T^nF)(a_1,\dots,a_n) = T(F(a_n,\dots a_1))$$ is a cochain morphism.
Therefore, with the preceding notation an antimorphism $\phi: B\to
A$ induces a morphism of cohomology $H^*(A, M) \to H^*(B^{\op},M)
\cong H^*(B, M^{\op})$ and an antimorphism $T:M\to N^{\op}$ induces
a morphism $H^*(A,M)$ \mbox{$\to H^*(A^{\op}, N^{\op})$}$\cong
H^*(A,N)$.

The next section uses the fact that the dual, $M^{\vee}= \operatorname{Hom}_k(M,k)$ of an $A$-module
$M$, is an $A^{\op}$-module in a natural way, and the
foregoing will imply that $H^*(A^{\op}, M^{\vee}) \cong H^*(A,
M^{\vee\,\op})$. However, even when $A$ is commutative and $M$ a
symmetric bimodule, $H^*(A, M)$ is generally \emph{not} isomorphic
to $H^*(A, M^{\vee})$, as an example there will show.

The dual of a Lie algebra $\g$ is similarly defined by $\g^{\vee} =
\operatorname{Hom}(\g, k)$ but in this case we may simply view
$\g^{\vee}$ as a $\g$-module: if $c \in \g$ and $f\in \g^{\vee}$
then $[c,f]$ is defined by setting $[c,f](a) = -f([c,a])$ for all $a
\in \g$.\, (If $V, W$ are modules over a Lie algebra $\g$ then
$\operatorname{Hom}_k(V,W)$ becomes a $\g$-module as follows: if
$\phi:V\to W$ and $c\in\g$ then $[c,\phi]$ is defined by setting
$[c,\phi](a) = [c,\phi a] - \phi([c,a])$ for all $a\in \g$. The
preceding is the special case where $V=\g, W= k$; note that the
coefficient ring $k$ is a Lie module over any Lie algebra $\g$ but
the operation is trivial, so the first term on the right, which
would be $[c, f(a)]$, vanishes.) Theorem 1 is not really necessary
in the Lie case where cocycles are alternating: Sending a cocycle to
its opposite leaves it unchanged if the dimension is even and just
reverses the sign if it is odd.

\section{Self-dual and Frobenius algebras}
Suppose that $A$ is an associative $k$-algebra and let $M$ be an
$A$-bimodule which (like $A$) will always be assumed to be free and
of finite rank as a $k$-module. Its \emph{dual}, $M^{\vee} =
\operatorname{Hom}_k(M, k)$ is then again a free $k$-module and of
the same rank as $M$ but should be viewed as an $A^{\op}$-bimodule
since the action of $A$ is reversed: If $f\in M^{\vee}$ and $a,b,x
\in A$ then $(afb)(x) = f(bxa)$. Thus $M^{\vee\,\op}$ is
again an $A$-bimodule. Note that an $A$-bimodule morphism $M\to N$
induces an $A^{\op}$-module morphism $N^{\vee} \to M^{\vee}$. A
\emph{self-dual} $A$-bimodule $M$ is one with an $A$-bimodule
isomorphism $\rho:M \to M^{\vee\op}$ or equivalently an
antiisomorphism $\rho^{\op}:M \to M^{\vee}$. The latter then gives
rise to a non-degenerate bilinear form $<-,->:M\times M \to k$ by
setting $<m,m'> = (\rho\,m)(m')$. However, while the existence of a
non-degenerate form $<-,->:M\times M \to k$ gives rise  to a 
$k$-module monomorphism $\rho:M\to M^{\vee}$ by sending $m$ to $<m,->$,
it is a bimodule morphism if and only if $<amb, m'> = <m, bm'a>$ for
all $a,b \in A$.  In general this induced $\rho$ need not be an
epimorphism unless $k$ is a field and $M$ is finite-dimensional. For
simplicity we always assume that $M$ is free and of
finite rank over the coefficient ring $k$. A sufficient condition
that $\rho$ be onto is then that there exist a basis $\{m_1,\dots,m_k\}$
such that $\det<m_i,m_j>$ is invertible in $k$, or equivalently,
that there exist dual bases $\{m_1,\dots,m_k\}$ and
$\{m_1',\dots,m_k'\}$ for $M$, i.e., ones with $<m_i,m_j'> =
\delta_{ij}$. (The problem when $k$ is not a field can be
illustrated by taking $k = \mathbb{Z} = M$. With the bilinear form
$<m,n>= mn$ the module is actually self dual; $\{1\}$ is a basis and
is self dual. However, with the form $<m,n> = 2mn$ there no longer
exist dual bases -- one would like to take $1^{\vee} = 1/2$ but that
is not in $\mathbb{Z}$.)

Viewing $A$ as a bimodule over itself, we may in particular consider
$A^{\vee\,\op}$ and call $A$ \emph{quasi self-dual} if there is a
canonical isomorphism of graded cohomology modules $H^*(A,A) \cong
H^*(A, A^{\vee\op})$. If we have a morphism $B \to A$ of such
algebras then the sequence
\begin{multline*}
H^*(A,A)\cong H^*(A, A^{\vee\,\op}) \to H^*(B,A^{\vee\,\op}) \cong
H^*(B^{\op},A^{\vee})\\
\to H^*(B^{\op},B^{\vee})\cong H^*(B, B^{\vee\op}) \cong H^*(B,B)
\end{multline*}
 exhibits $H^*(A,A)$ as a contravariant functor of $A$. This
principle is not restricted to associative algebras but clearly
holds equally well, e.g., for Lie algebras. 

The direct sum of quasi
self-dual algebras is again such. Since we are considering only
algebras of finite rank (partly to avoid topological problems) we
have $(A_1 \otimes A_2)^{\vee} = A_1^{\vee} \otimes A_2^{\vee}$.
This, together with the fact that in general if $M_1, M_2$ are $A_1,
A_2$ bimodules, respectively, then $H^*(A_1\otimes A_2, M_1 \otimes
M_2) = H^*(A_1, M_1) \otimes H^*(A_2,M_2)$ shows that the category
of quasi self-dual algebras is closed under tensor products. It is
clearly also closed under direct sums but is not closed under taking
quotients, as will be seen.

Suppose now that we have a $k$-module morphism $\phi:A\to A^{\vee}$
which for the moment need be neither a monomorphism nor epimorphism.
Then we can define a bilinear form $<-,->:A\times A \to A$ by $<a,b>
= (\phi a)(b)$, and conversely. The condition that $\phi$ be an
antimorphism from $A$ viewed as an $A$-bimodule to $A^{\op}$ viewed
as an $A^{\op}$-bimodule then is equivalent to having both
\begin{equation}<ac,b> =<a,cb> \quad \text{and} \quad <ca,b> = <a,bc>.\label{eq:dual}\end{equation}
If the algebra is unital (which we generally assume) then the form
must be symmetric for we have
$$<c,b> = <1\cdot c, b> = <1,
cb> = <b\cdot1,c> = <bc,1> = <b,c>.$$ The conditions (\ref{eq:dual})
are equivalent if the form is symmetric, and that in turn will be
the case whenever the algebra is commutative, but a non-commutative
algebra may still have a symmetric form, e.g.  finite group rings
(below). If a bilinear form for which the associated $k$-linear
mapping $A \to A^{\vee}$ is an isomorphism satisfies the conditions
(\ref{eq:dual}) then that form will be called \emph{dualizing} and
an algebra with a dualizing form will be called \emph{self-dual}
since these are precisely the ones with an $A$-bimodule isomorphism
$A \to A^{\vee\op}$; they are in particular quasi self-dual. When
$A$ is graded then by ``commutative'' we will always mean
commutative in the graded sense (sometimes called
``supercommutative''). In that case, if $\deg a = r,\, \deg b = s$
then a `symmetric' form must have $<a,b> = (-1)^{rs}<b,a>$.

The classical definition of a \emph{Frobenius algebra} $A$ (cf.
Nakayama, \cite{Nak:1939}) is one which is a finite dimensional
associative algebra over a field $k$ together with a linear
functional $f:A\to k$ satisfying any of the following conditions
which are equivalent in the presence of finite dimensionality: (i)
the kernel of $f$ contains no left ideal, (ii) the kernel of $f$
contains no right ideal, (iii) the bilinear ``Frobenius form''
$<-,->:A\times A \to k$ defined by $<a,b> = f(ab)$ is
non-degenerate. The form $<-,->$ then has the property that
\mbox{$<a,bc>=<ab,c>$}, and conversely, the existence of a
non-degenerate bilinear form with this property implies that $A$ is
Frobenius by setting $f(a) = <a,1>$.  A \emph{symmetric Frobenius}
algebra is one for which the form is symmetric.

\begin{theorem} A finite dimensional associative algebra over a
field is self-dual if and only if it is a symmetric Frobenius
algebra. In particular commutative Frobenius algebras are
self-dual.
\end{theorem}
\textsc{Proof.} With the above notation we now have both $<a,cb> =
<ac,b>$ (automatic from the definition) and $<a,bc> =<bc,a> =
<b,ca>=<ca,b>$, using the symmetry, so the Frobenius form is
dualizing. $ \Box$
\medskip

A sufficient condition for symmetry of a Frobenius algebra $A$ is
that it possess an involution $\sigma$ (i.e., antiautomorphism whose
square is the identity) preserving the Frobenius form i.e., with
$<\sigma a, \sigma b> = <a,b>$ (or equivalently, if the form is
defined by the functional $f:A\to k$, with $f(\sigma a) = f(a)$).
For then we have
\begin{multline*}
<a, bc> = <\sigma a, \sigma(bc)> = <\sigma a, \sigma
c\cdot\sigma b>\\
 = <\sigma a\cdot\sigma c, \sigma b> = <\sigma(ca),
\sigma b> = <ca,b>.
\end{multline*}

There are important examples of symmetric Frobenius algebras. The
group algebra $kG$ of a finite group $G$ is Frobenius: define
$<a,b>$ to be the coefficient of the identity $1=1_G$ in the product
$ab$. This, however, is the same as the coefficient of $1$ in $ba$,
so $kG$ is symmetric. In fact, here we do not have to assume that
$k$ is field, for if we take as a basis for $kG$ the elements of $G$
then the dual basis consists of their inverses. Thus the group ring
$kG$ of a finite group $G$ over any (commutative, associate, unital)
ring $k$ is self-dual. The map sending $g$ to $g^{-1}$, extended
linearly, is an involution preserving the form. Another example is
the de Rham cohomology ring of a compact manifold $\mathcal{M}$
where, if $a, b$ are cocycles then one sets
$$<a,b> = \int_\mathcal{M}a\wedge b\quad .$$ We should like to be
able to do the same for cohomology with integer coefficients,
defining $<a,b>$ to be the evaluation of $a\smile b$ on the
fundamental cycle, but in addition to the general problem when the
coefficient ring is not a field there may now be torsion. In some
favorable cases, however, this is still possible.

Commutative Frobenius algebras have recently been shown to play an
important role in the algebraic treatment and axiomatic foundation
of topological quantum field theory since such an algebra determines
uniquely (up to isomorphism) a 1+1 dimensional TQFT. More
precisely, the category of commutative Frobenius $k$ algebras is
equivalent to the category of symmetric strong monoidal functors
from the category of 2-dimensional cobordisms  between 1-dimensional
manifolds to the category of vector spaces over $k$, \cite{Abrams:Frobenius}.

Closely related to the group algebras of finite commutative groups
are the algebras of the form $k[t]/t^{n+1}$ and their tensor
products. For let $C_q$ denote the (multiplicative) cyclic group of
order $q$, let $a$ be a generator,  $p$ be a prime, and set $q=p^r$.
If $k$ has characteristic $p$, then $kC_{q} \cong k[t]/t^q$, the
isomorphism being given by $1-a \mapsto t$. (Since for finite groups
$G_1$ and $G_2$ we have $k(G_1 \times G_2) = kG_1 \otimes kG_2$,
when dealing with commutative groups one need only consider the
cyclic case.)

Algebras of the form $k[t]/t^{n+1}$ are Frobenius. For letting
$u_0,\dots,u_n$ be the dual basis to $1, t,\dots, t^n$ we have
$t\cdot u_i = u_{i-1}, i = 1,\dots,n;\,\, t\cdot u_0 = 0$, so we can
set $<t^i, t^j> = 1$ if $i+j=n,$ and set it equal to $0$ otherwise.
It follows, e.g., that $k[x]/x^r \otimes k[x]/x^s$ is Frobenius for
all $r$ and $s$, hence self-dual. Taking $r = s = 2$ it follows that
$k[x,y]/(x^2, y^2)$ is Frobenius. This is free of rank 4 over $k$,
therefore of dimension 4 when $k$ is a field. However, its
3-dimensional quotient $A = k[x,y]/(x^2, y^2, xy)$ is not self-dual,
for one can check that $H^1(A,A)$ has dimension 4 but $H^1(A,
A^{\vee})$ has dimension 3. (Note here that $A= A^{\op}$ and that
$A^{\vee}$ is symmetric.) It follows that $H^*(A, A^{\vee\op})$ is
not isomorphic to $H^*(A, A)$ so $A$ is not quasi self-dual, and in
particular not Frobenius. This example shows that a quotient of a
quasi self-dual algebra need not be quasi self-dual.

That $k[t]/t^{n+1}$ is Frobenius may be viewed as one of the
``favorable cases '' mentioned above: The cohomology of complex
projective $n$-space $\mathbb{CP}^n$ is isomorphic to $k[t]/t^{n+1}$
where $t$ is the cohomology class of dimension 2 defined by any
hyperplane. (Since the dimension is even the algebra is actually
commutative in the ordinary sense.) In this case the form mentioned
above, evaluation of the cup product on the fundamental cycle, is
indeed dualizing.

If $G$ is a group whose order is invertible in $k$ then $H^n(kG, M)
= 0$ for all $n > 0$ and all $kG$ modules $M$ (a form of Maschke's
theorem, cf. the next section). However, for $k = \Z$ these cohomology
groups are generally not trivial although they are all necessarily
torsion modules. (Computing the cohomology groups $H^n(\Z C_n, \Z)$
is simplified by the fact that that $\Z$, viewed as a
trivial $\Z C_n$ module, has a projective
resolution which is periodic of order 2 and in which every
projective module is just $\Z C_n$ itself. This property is shared
by all $k[t]/t^{n+1}$. We have for all $k$
\begin{equation*}
\begin{CD}\dots
k[t]/t^{n+1}@>\partial_i>>k[t]/t^{n+1}@>\partial_{i-1}>>\dots@>
\partial_1>>k[t]/t^{n+1}@>{\epsilon}>>k
\end{CD}
\end{equation*}
where $\epsilon$ is reduction mod $t$ and $\partial_i$ is
multiplication by t for $i$ odd and multiplication by $t^n$ for $i$
even.)

The graded Lie structure on $H^*(A,A)$ when $A = kG$ or $A = k[t]/^{n+1}$ 
is generally not trivial as one can see already by considering the
commutator of derivations.  In the case of $k[t]/t^{n+1}$, for
example, if $D_1t=t^r$ and $D_2t=t^s$ then $D_1D_2 = st_{r+s-1}$, so
$[D_1,D_2]t = (s-r)t^{r+s-1}$, which is generally not 0. (We must
assume that $1 \le r,s \le n$ except, e.g. in characteristic $p$
where for example if $A = k[t]/t^{p^r}$ one can allow $Dt = 1$.)
Both $H^*(kG,kG)$ and $H^*(kG,k)$ are rings with multiplications
induced by those in the respective coefficient modules $kG$ and $k$,
and the module morphism $kG \to k$ is in fact a ring morphism; it
induces a ring morphism $H^*(kG,kG) \to H^*(kG,k)$, but the graded
Lie part of the Gerstenhaber algebra structure on $H^*(kG,kG)$ is
lost.

\section{The Lie case} The basic result for Lie algebras is the
following
\begin{theorem} Lett $\g$ be a Lie algebra with a bilinear form $<-,->$ for
which the associated mapping $\phi: a\to a^{\vee}$ is a $k$ linear
isomorphism. Then the form is  dualizing if and only if it is
invariant, i.e., if and only if $<[c,a], b> + <a,[c,b]> = 0$ for all $a,b,c
\in \g$.
\end{theorem}
\textsc{Proof.} Suppose that $\phi:\g \to \g^{\vee}$ is a $k$-module
isomorphism giving rise to the form $<a,b> = (\phi a)(b)$. To be a
module morphism we must have $(\phi([c,a])(b) = [c,\phi a](b)$ for
all $a,b,c$. The left side is $<[c,a],b>$; the right (recalling that
$\phi a$ is a $k$-module map from $\g$ to the trivial $k$-module) is
$-<a, [c,b]>$. $\square$
\medskip

In the Lie case a  dualizing form need not be symmetric, as we shall
see by example. However, we have

\begin{lemma} If a  Lie algebra $\g$ of characteristic $\ne 2$ has a
skew dualizing form then $\g$ is Abelian.
\end{lemma}
\textsc{Proof.} Suppose that $<-,->$ is a skew dualizing form on
$\g$. Then$\mbox{$<[c,b], a>$} = -<b,[c,a]> = <[c,a],b>$, but the
left term equals $<c,[b,a]>$ and the right equals $<c,[a,b]>$.
Assuming that the characteristic is not $2$ this implies that
$<c,[a,b]> = 0$ for all $a,b,c$\,; since the form is non-degenerate
one has $[a,b] = 0$ for all $a$ and $b$, so the Lie algebra is
Abelian. $\square$

The close relationship between associative self-dual and Frobenius
algebras no longer holds in the Lie case. A \emph{Frobenius Lie
algebra} $\g$ is one which is finite dimensional over a field and
has a functional $f:L\to k$ such that the skew bilinear
\emph{Kirillov form}  $<a,b> =f([a,b])$ is non-degenerate. The
Kirillov form is by definition a coboundary in the
Chevalley-Eilenberg theory; a quasi-Frobenius Lie algebra is one
with a skew non-degenerate `Kirillov' form which is just a
$2$-cocycle. By the foregoing the Kirillov form can not be
dualizing. Conceivably a Frobenius or quasi-Frobenius Lie algebra
could still be self-dual, but relative to some form other than the
Kirillov form.

The Killing form of a finite dimensional semisimple real or complex
Lie algebra is non-degenerate and invariant so these algebras are
self-dual. However, the cohomology of any such algebra with
coefficients in any finite dimensional module other than the trivial
module vanishes. For an example with non-trivial cohomology, let $V$
be a real vector space with an inner product $<-,->_V$, $\G$ be its
orthogonal group and $\g$ be the Lie algebra of $G$. Then $V$ is a
$\g$ module and $<-,->_V$ is invariant. Now consider the semi-direct
product $\g \ltimes V$  (a split extension of $\g$ by $V$); its Lie
multiplication is given by $[(a,v), (b,w)] = ([a,b], [a,w] -
[b,v])$. The cohomology of such a Lie algebra with coefficients in
itself is generally not zero cf.  \cite{Richardson:rigidity}.
Denoting the Killing form on $\g$ by $<-,->_K$, the symmetric form
on $\g \ltimes V$ defined by $<(a,v),(b,w)> = \mbox{$<a,b>_K$}
+<v,w>_V$ is non-degenerate and invariant, hence dualizing. For an
example of a case where the dualizing form is neither symmetric nor
skew, let $V$ be of even dimension with a non-degenerate skew form,
$\G$ be its symplectic group, and perform the same construction. 
There may be other examples of self-dual Lie algebras, but note
 that a finite dimensional real or complex Lie
algebra with a non-degenerate invariant form and no one-dimensional
ideal is in fact just a semi-direct product of a semisimple Lie
algebra with some module over that algebra having no irreducible
component of dimension one. This is a consequence of \cite[Theorem
3, p.71]{Jacobson:LieAlgs} and the fact that the cohomology of a
semisimple Lie algebra with coefficients in a module whose
decomposition into simple components has no component of dimension
$1$ vanishes in every positive dimension, and in particular, in
dimension $2$.

\section{Separable algebras and relative Hochschild cohomology}
A unital $k$-algebra $S$  is \emph{separable} (over $k$) if it is
projective in the category of $S$-bimodules. Viewing both $S$ and $S
\otimes S^{\op}$ as $S$-bimodules, the multiplication map $\mu: S
\otimes S^{\op} \to S$, which is an $S$-bimodule morphism, then has
a splitting $\nu: S \to S\otimes S^{\op}$, i.e., an $S$-bimodule
morphism such that $\mu\nu\mu = \mu$. Since $S\otimes S^{\op}$ is a
free bimodule of rank 1 over $S$, the existence of such a splitting
is equivalent to separability, since it exhibits $S$ as a direct
summand of a free module (for note that $S\otimes S^{\op}$ is a free
$S$-bimodule  with generator $1\otimes1$). 
 If we have such a bimodule splitting, for the
moment write $\nu(1_S) = \sum\,x_i\otimes y_i = e_{\text{sep}}$ or
simply $e$. Since $s\nu(1) = \nu(s\cdot 1) = \nu(1\cdot s) =
\nu(1)s$ for all $s\in S$, this $e$ has the remarkable property that
$$\sum\,sx_i\otimes y_i= \sum\,x_i\otimes y_is, \qquad \mathrm{all} \quad s\in S.$$
Since $1 = \mu(e)= \sum\,x_iy_i$ it follows that $e$ is idempotent;
it is called a \emph{separability idempotent} for $S$. (These are generally not unique.) 
Using it we
can show that any morphism $f:M\to N$ of left $S$-modules which
\emph{a priori} splits only as a morphism of $k$-modules actually
splits as a morphism of left $S$-modules. For if $g:N\to M$ is a 
$k$-morphism such that $fgf = f$ then $\bar g:N\to M$ defined by $\bar
g(n) = \sum_i x_ig(y_in)$ is a left $S$-module morphism such that
$f\bar g f= f$. If $g$ is already a left module morphism then $\bar
g$ will be identical with $g$, so $e$ projects $k$-module splittings
onto $S$-module splittings. The same holds for right modules and
also for bimodules by similar arguments. The last is equivalent to
separability for $\mu$ always has a $k$-splitting: One can send $s\in S$ to
$s\otimes 1$, and this can by hypothesis then be projected onto a
bimodule splitting. The one-sided conditions, however, are not
strong enough, since if $k$ is a field and $S$ an inseparable
extension then a left $S$-module is just a vector space, so the left
and right splitting properties both hold. But $S\otimes S$, will
then contain a non-trivial radical and the bimodule splitting
property will not hold. The term ``separable'' derives ultimately
from the fact that a finite field extension is separable in the
classical sense precisely when it is so in that above. (The present
definition is due to M. Auslander and O. Goldman,
\cite{AusGold:Brauer}, based on remarkable previous work by Azumaya,
\cite{Azumaya}.)

Using the definition of cohomology by projective resolutions, it
follows immediately from the two-sided splitting condition that if
$S$ is separable over $k$ then $H^n(S,M) = 0$ for all $n\ge 1$ and
all $S$-bimodules $M$. The latter is, in fact, another equivalent
criterion for separability. By dimension shifting techniques one
need in fact only assume that $H^1(A,M) = 0$ \,for all $M$ and it is
this latter that is equivalent to the two-sided splitting condition.
The vanishing of cohomology is the basic property of separability
that we will need.

It is usually easiest to prove separability by exhibiting a
separability idempotent. The algebra of $n \times n$ matrices over
$k$ is separable; its separability idempotent is
$\Sigma_i\,e_{i,1}\otimes e_{1,i}$. (The fixed internal index $1$
can be replaced by any other.) The group ring $kG$ of a finite group is
separable over $k$ if the order of $G$ is invertible in $k$; if the
order is $N$ then $e = \frac{1}{N} \sum_{g\in G}g\otimes g^{-1}$.
(This, in essence, is Maschke's theorem.) Tensor products and finite
direct sums of separable algebras are separable. In particular, $k$
is separable over itself and therefore a finite direct sum of copies
of $k$ is separable, a fact which we will need. A separable 
$k$-algebra which is projective as a module over $k$ must be finitely
generated; in particular, an algebra which is separable over a field
is finite dimensional. (For a readable discussion of separability,
cf. \cite{DeMeyerIng:SepAlgs}.)

The use of separability can sometimes simplify the computation of
Hochschild cohomology. First, call an $n$-cochain $F$ in $C^n(A,M)$
\emph{normalized} if it vanishes whenever any of its arguments is
the unit element of $k$. These cochains form a subcomplex of the
full Hochschild complex whose inclusion into the full Hochschild
complex induces an isomorphism of cohomology. Suppose now that $S$
is a $k$-subalgebra of $A$, arbitrary except that we will always
assume that the unit element of $A$ is contained in $S$. An
$S$-\emph{relative} cochain $F\in C^n(A,M)$ is one such that for all
$a_1, \dots, a_n \in A$ and $s\in S$ we have
%\begin{equation*}
\begin{align*}
F(a_1,\dots, a_is,a_{i+1},\dots,a_n) &= F(a_1, \dots, a_i,sa_{i+1},\dots,a_n), &\quad i =1, \dots,n-1,\\
F(sa_1,\dots,a_n) &= sF(a_1,\dots,a_n),&\text{and}\quad\quad  \\
F(a_1,\dots,a_ns) &= F(a_1,\dots,a_n)s.&
\end{align*}
%\end{equation*}
If, moreover,  $F$ is normalized then it must vanish whenever any
argument is in $S$.  The relative cochain groups, denoted $C^n(A, S;
M)$, also form a subcomplex of the Hochschild complex. The result we
need is that when $S$ is a separable algebra over $k$ the inclusion
of the complex of  $S$-relative cochains into the full Hochschild
complex induces an isomorphism of cohomology. This may be difficult
to see from the original Hochschild definition of cohomology but is
relatively transparent from the point of view of projective
resolutions since separability actually allows one to take tensor
products over $S$ rather than $k$ in the bar resolution. Finally,
the normalized relative cochain groups $\overline C^n(A,S;M)$ form a
subcomplex of the relative groups and their inclusion into the full
Hochschild cochain complex again induces an isomorphism of
cohomology.  It is this last subcomplex of normalized $S$-relative
cochains which will be essential in the next section.

\section{Poset algebras}
Let $\I = \{\dots, i, j, \dots\}$ be a finite poset of order $N$
with no cycles and partial order denoted by $\pr$. Extending the
partial order to a total order we may assume that $\I
=\{1,\dots,N\}$ where the partial order $\pr$ is compatible with the
natural order. Associated to $\I$ is the algebra $A = A(\I)$ of
upper triangular matrices spanned by the matrix units $e_{ij}$ with
$i\pr j$; these are closed under multiplication. The subalgebra $S$
spanned by the $e_{ii}, i = 1,\dots, N$ has the same unit as $A$ and
is isomorphic to a direct sum of copies of $k$, hence is separable,
so we may compute the Hochschild cohomology of $A$ with coefficients
in any module, and in particular  $H^*(A,A)$ and $H^*(A,A^{\op})$,
using $S$-relative cochains.

Viewing $\I$ as a small category in which $\operatorname{Hom}(i,j)$
consists of a single morphism $i \to j$ if $i\pr j$ (the identity
morphism when $i = j$) and is empty otherwise, its \emph{nerve}
$\Sigma = \Sigma(\I)$ is the simplicial complex whose $n$-simplices
$\sigma$ are the $n$-tuples of composable morphisms $i_0\to i_1\to
\dots\to i_n$,  whose module $C_n$ of $n$-chains consists of the
linear combinations of these, and where $\partial_0\sigma$ just
omits the first morphism, $\partial_n\sigma$ omits the last and
$\partial_r\sigma, \, 0<r<n$ is obtained by replacing $i_{r-1}\to
i_r\to i_{r+1}$ by $i_{r-1} \to i_{r+1}$. (The $0$-simplices or
vertices are just the elements of $\I$ and $\partial_0(i_0\to i_1) =
i_1,\,\partial_1(i_0\to i_1) = i_0$.) One sets $\partial\sigma =
\Sigma_{i=0}^n(-1)^r\partial_r\sigma$. Denote the simplicial
cohomology of $\Sigma =\Sigma(\I)$ with coefficients in $k$ by
$H^*(\Sigma, k)$. With these notations we can recapitulate the
theorem which is the basis of the more general result in
\cite{GS:Simplicial}.
\begin{theorem} There is a canonical isomorphism $H^*(A,A) \cong
H^*(\Sigma, k)$.
\end{theorem}
{\textsc Proof.} We can compute the left side using $S$-relative
cochains $F \in C^n(A, S;A)$.  A cochain is completely determined
when its arguments are amongst the $e_{ij}$. If $F$ is $S$-relative
then
\begin{equation*}
F(\dots, e_{ij}, e_{kl},\dots) = F(\dots, e_{ij}e_{jj},
e_{kl},\dots)= F(\dots, e_{ij}, e_{jj}e_{kl}, \dots).
\end{equation*}
This vanishes if $j\ne k$ so the only non-zero values of $F$ are
those of the form $F(e_{i_0i_1}, e_{i_1i_2}, \dots,
e_{i_{n-1}i_n})$, where $i_0 \preceq i_1 \preceq \cdots \preceq
i_n$. Also, since $F$ is $S$-relative these values must lie in
$e_{i_0i_0}Ae_{i_ni_n}$. The latter is a free module of rank 1
spanned by $e_{i_0i_n}$ so the value is $\lambda e_{i_0i_n}$ for
some $\lambda \in k$. Thus  $F$ assigns to every $\sigma = (i_0\to
i_1\to\dots \to i_n)$ an element $\lambda \in k$ and so may be
viewed as an element of $C^n(\Sigma, k)$. Conversely, every such
simplicial cochain defines a unique $S$-relative Hochschild cochain.
It is easy to check that the simplicial and Hochschild coboundaries
then correspond. $\square$
\medskip

Any simplicial complex gives rise to a poset whose objects are its
faces with the relation ``face of''. The nerve of this poset is just
the barycentric subdivision of the the original simplicial complex
and hence has cohomology isomorphic to that of the original.

\begin{corollary} For every finite simplicial complex $\Sigma$ and
coefficient ring $k$ there is a poset $k$-algebra $A$ such that
$H^*(A,A) \cong H^*(\Sigma, k).$ $\square$
\end{corollary}

It is not difficult to see that the isomorphism respects the cup
product, so the fact that the cup product in simplicial cohomology
is graded commutative is actually a consequence of the more general
fact that it is so in $H^*(A,A)$ for every associative algebra $A$,
as first shown in \cite{G:Cohom63}.

Let$ \{e_{ij}^*\}$ denote the dual basis to $\{e_{ij}\}$. For the
operation of $A^{\op}$ on the dual bimodule $A^{\vee}$ a simple
computation shows that
\begin{align*}
e_{ik}^*e_{lj} & = 0 &\text{ if } \quad &i\ne l \text{ or }  j
\not\pr k;\quad
e_{ik}^*e_{jk} & = e_{jk}^* \quad &\text{ if } \quad &i\pr j\pr k\\
e_{ij}e_{hl}^* & = 0 &\text{ if } \quad &j\ne l \text{ or } h
\not\pr i; \quad e_{ij}e_{hj}^* & = e_{hk}^*\quad  &\text{ if }
\quad &h\pr i \pr j
\end{align*}

Here for the moment we write  simply $e_{ik}^*e_{lj}$ (and not
$e_{ik}^*\circ e_{jk}$) since in the foregoing we have just the
natural operation of $A^{\op}$ on $A^{\vee}$.

\begin{theorem}
If $A$ is a poset algebra then there is a natural isomorphism
$$H^*(A,A) \cong H^*(A, A^{\vee\,\op});$$ poset algebras are quasi
self-dual.
\end{theorem}
{\textsc Proof.} Using $S$-relative cochains as before we find that
$F(e_{i_0i_1}, e_{i_1i_2}, \dots, e_{i_{n-1}i_n})$ lies in
$e_{i_0i_0}\circ A^{\vee\,\op}\circ e_{i_ni_n}=
e_{i_ni_n}A^{\vee}e_{i_0i_0}$ which is again a free module of rank
1, spanned by $e_{i_0i_n}^*$. The rest follows as before. $\square$
\medskip

Since poset algebras over a fixed ring $k$ are quasi self-dual any
algebra morphism between them (or between one and any quasi
self-dual algebra) induces a morphism of cohomology. Not all
morphisms between poset algebras need be induced by morphisms
between the underlying posets but we will not examine here which
others are possible. As mentioned, they are probably very
restricted.

What has been proven here for poset algebras actually applies to a
wider class of algebras called ``triangular'' in
\cite{GS:Triangular}; these are deformations of poset algebras.
However, from the deformed algebra one can reconstruct the original
poset, and hence the original poset algebra. The cohomologies of the
deformed algebra and of the original algebra (each with coefficients
in itself) are identical and these are identical with the simplicial
cohomology of the nerve of the poset. Thus in the very special case
of poset algebras deformation does not alter cohomology. It follows
that triangular algebras in the sense of \cite{GS:Triangular} are
also quasi self-dual.  We may view any space
with a finite triangulation as having a hidden algebra structure
depending on which deformation of its poset algebra we associate to
it, but the family of these is already completely determined by the
cohomology of the space alone. Perhaps for these spaces the concept
of space should mean not just the topological object but a pair
consisting of that object together with the algebra associated to
it.

Classifying self-dual algebras seems difficult since this includes 
classifying 1+1 dimensional topological quantum field theories (but see the cited works
 of Abrams and Sawin). Any classification would have to be into families since
 Frobenius alebras can deform in non-trivial ways; a simple example is given in 
the next section. Classification might be aided if one could determine when a quasi self-dual
algebra is related in some way to a geometric object, e.g.
 as with a poset algebra and the nerve of its poset, or even a
more extended way  as in the relation between commutative Frobenius
algebras and topological quantum field theories.

\section{Deformation of $\C[t]/t^{n+1}$ as a Frobenius algebra. }
A natural question concerning any algebraic structure is, What is 
its deformation theory? In particular,
is a deformation of a self-dual algebra again self-dual and a
deformation of a Frobenius algebra again Frobenius? This seems
to be a difficult question in general but one can give a positive answer 
at least for certain deformations of the commutative Frobenius 
algebra $A = \C[t]/t^{n+1}$
encountered earlier. Since this algebra is generated over $\C$ by  a single
element the same must be true of any deformation. Any deformation must therefore be
equivalent to one of the form $A_{\hbar} = \C[t]/(t^{n+1} - \hbar
p(t))$ where $p(t)$ is a polynomial of degree $\le n$ and $\hbar$ is
the deformation parameter. (A quotient of a fixed algebra by a
varying ideal depending on some parameter $\hbar$ can not always be
exhibited as a deformation in the sense of \cite{G:DefI} since the
dimension may change but there is no problem here since one has
a fixed basis for the quotient algebra, namely $\{1, t, t^2, \dots,
t^n\}$, independent of the value of $\hbar$.) 

In general, the coefficients of the
polynomial $p(t)$ will themselves be power series in $\hbar$. 
For simplicity we consider here only the case where the deformation is defined
by a polynomial $p(t)$ with constant coefficients and ask the following question:
Can the Frobenius form on $\C[t]/t^{n+1}$
given earlier, namely that in which $<t^i,t^j> = 1$ if $i+j = n$ and is $0$
otherwise, be deformed along with the algebra so that the deformation is 
again a commutative Frobenius algebra? Denoting the deformed alebra by
 $A_h$, if it does carry a deformed form then  the values of
 $<1, t^i>,\, i = 0, \dots, n$ completely
determine the form. For simplicity, write $<t^i,t^j> = c_{i,j}$ and
set $p(t) = a_0 + a_1t+ \dots a_nt^n$. Note that
\begin{equation*}
<t^i, t^j> = <t\cdot t^{i-1}, t^j> = <t^{i-1},t\cdot t^j> =
<t^{i-1}, t^{j+1}>
\end{equation*}
so the value of $c_{i,j}$ depends only on $i+j$.  Then from $<t^i,
t^n> = <t^{i-1}, t^{n+1}> = \,\,<t^{i-1}, \hbar p(t)>$ we get the
recursion
\begin{equation*}
c_{i,n} = \hbar(c_{i-1,0}a_0 + c_{i-1,1}a_1 + \dots, c_{i-1,n}a_n).
\end{equation*}
The values of $c_{0,0}, c_{0,1}, \dots, c_{0,n-1},  c_{0,n}$
therefore determine all the $c_{i,j}$. Setting these equal to $0,0,
\dots, 0,1$, respectively, the matrix $||c_{i,j}||$ of the resulting
form has $1$ everywhere on the main antidiagonal, zero entries above
the main antidiagonal and polynomials in $\hbar$ without constant
term below the main antidiagonal. Its determinant is therefore
$(-1)^{\lfloor \frac{n+1}{2}\rfloor}$, so the resulting form is
non-degenerate (and would be with any coefficient ring). The form on
$A_{\hbar}$ varies continuously with $\hbar$ and reduces to the
original simple form on $A$ when $\hbar = 0$, so in this simple case 
the algebra
deformation from $A$ to $A_{\hbar}$ has indeed `dragged along' a
deformation of the Frobenius form on $A$ to one on $A_{\hbar}$.

While the foregoing shows that $\C[t]/t^{n+1}$ with the given form
indeed deforms as a Frobenius algebra, the deformation exhibited actually proceeds by jumps.
For suppose that $P(t)$ is a monic polynomial in $\C[t]$ of degree
$n+1$ and that it factors into $P(t) =
(t-r_1)^{n_1}(t-r_2)^{n_2}\cdots (t-r_k)^{n_k}$. Then by the Chinese
Remainder Theorem,
\begin{equation*}
\begin{split}
\C[t]/P(T) \cong \C[t]/(t-r_1)^{n_1} \oplus \C[t]/(t-r_2)^{n_2}
\cdots \oplus \C[t]/(t-r_k)^{n_k}\\ \cong \C[t]/t^{n_1} \oplus
\C[t]/t^{n_2} \oplus\cdots \oplus \C[t]/t^{n_k}.
\end{split}
\end{equation*}
The structure of $\C[t]/P(T)$ therefore depends only on the
partition $\{n_1, n_2, \dots, n_k\}$ of $n+1$, where we may assume
without loss of generality that
  $n_1 \ge n_2 \ge \dots \ge n_k$. Since for small non-zero values of
 $\hbar$ the partition associated with the factorization of $P(t) =
 t^{n+1} -\hbar p(t)$ is constant, the deformation $A_{\hbar}$ is in
 fact a jump deformation; all sufficiently small non-zero values of
 $\hbar$ give isomorphic algebras. There is a natural partial order
 to partitions of an integer $n+1$, with $\{n+1\}$ the largest and
 $\{1^{n+1}\}$ (meaning $1$ repeated $n+1$ times) the smallest, and
 this partial order determines a partial order on the possible deformations
 of $\C[t]/t^{n+1}$, which is itself the highest in the partial
 order. The lowest is any $\C[t]/(t^{n+1} - a)$ with $a$ any non-zero
 complex number; these are all isomorphic to a direct sum of $n+1$
 copies of $\C$.

 Referring back to the associated TQFTs, they deform, but the
 deformations are jumps; one can not `see'   any intermediate
 stages. This raises the questions of whether (i) there is any intrinsic
 way to parameterize 1+1 dimensional TQFTs without referring to
 their associated algebras, and whether (ii) unlike the 
 deformations of $\C[t]/t^{n+1}$ exhibited here, 
 there are any families which
 do not consist exclusively of jump deformations.

 \nocite{HG:NATO}

%\bibliographystyle{plain}
%\bibliography{newbib.bib}

\begin{thebibliography}{10}

\bibitem{Abrams:Frobenius}
Lowell Abrams.
\newblock {Two-dimensional quantum field theories and Frobenius algebras}.
\newblock {\em J. Knot Theory Ramifications}, 5:569--587, 1996.

\bibitem{AusGold:Brauer}
M.~Auslander and O.~Goldman.
\newblock {The Brauer group of a commutative ring}.
\newblock {\em J. Algebra}, 4:220--272, 1960.

\bibitem{Azumaya}
G.~Azumaya.
\newblock {On maximally central algebras}.
\newblock {\em Nagoya Math. J.}, 2:{119--150}, 1951.

\bibitem{DeMeyerIng:SepAlgs}
Frank DeMeyer and Edward Ingraham.
\newblock {\em Separable Algebras over Commutative Rings}.
\newblock Number 181 in Lecture Notes in Mathematics. Springer-Verlag, {New
  York \ Heidelberg \ Berlin}, 1971.

\bibitem{G:Cohom63}
M.~Gerstenhaber.
\newblock {The cohomology structure of an associative ring}.
\newblock {\em Ann. of Math.}, 78:267--288, 1963.

\bibitem{G:DefI}
M.~Gerstenhaber.
\newblock {On the deformation of rings and algebras }.
\newblock {\em Ann. of Math.}, 79:59--103, 1964.

\bibitem{GS:Triangular}
M.~Gerstenhaber and S.~D. Schack.
\newblock {\em Triangular algebras}, pages 447--498.
\newblock in \cite{HG:NATO}, 1988.

\bibitem{GS:Simplicial}
M.~Gerstenhaber and S.D. Schack.
\newblock {Simplicial cohomology is Hochschild cohomology}.
\newblock {\em J. Pure and Appl. Algebra}, 30:143--156, 1983.

\bibitem{HG:NATO}
Michiel Hazewinkel and Murray Gerstenhaber, editors.
\newblock {\em Deformation Theory of Algebras and Structures and Applications}.
\newblock NATO ASI Series C: Mathematical and Physical Sciences. Kluwer,
  Dordrecht \ Boston \ London, 1988.

\bibitem{Jacobson:LieAlgs}
Nathan Jacobson.
\newblock {\em Lie Algebras}.
\newblock Interscience Publishers division of John Wiley and Sons, New York \
  London, 1962.

\bibitem{Nak:1939}
T.~Nakayama.
\newblock {On Frobeniusean algebras I}.
\newblock {\em Ann. Math}, 40:611--633, 1940.

\bibitem{Richardson:rigidity}
{R. W. Richardson, Jr. }.
\newblock {On the rigidity of semi-direct products of Lie algebras}.
\newblock {\em Pacific J. Math.}, 22:339--344, 1967.

\bibitem{Sawin:DirectSum}
Stephen Sawin.
\newblock {Direct sum decompositions and indecomposable TQFTs}.
\newblock {\em J. Math. Phys.}, 36(12):6673--6680, 1995.

\bibitem{Sawin:LinksQGsTQFts}
Stephen Sawin.
\newblock {Links, quantum groups and TQFTs}.
\newblock {\em Bull. Amer. Math. Soc. (N.S.)}, 33(4):413--445, 1996.

\end{thebibliography}
\end{document}